\documentclass[12pt,twoside]{article}
\usepackage{amsmath,amsfonts,amssymb,a4,enumerate,epsfig,array,theorem,psfrag}

\newtheorem{theo}{Theorem}

\newtheorem{prop}{Proposition}[section]
\newtheorem{coro}[prop]{Corollary}

\newcommand{\Diff}{\operatorname{Diff}}
\newcommand{\cone}{\bowtie}
\newcommand{\arccosh}{\operatorname{arc cosh}}
\newcommand{\sgn}{\operatorname{sgn}}

\newcommand{\tr}{\operatorname{tr}}

\newcommand{\NN}{{\mathbb{N}}}
\newcommand{\ZZ}{{\mathbb{Z}}}
\newcommand{\RR}{{\mathbb{R}}}

\newcommand{\Ss}{{\mathbb{S}}}

\newcommand{\HH}{{\mathbb{H}}}
\newcommand{\Cc}{{\cal{C}}}
\newcommand{\Pp}{{\cal{P}}}
\newcommand{\Ff}{{\cal{F}}}
\newcommand{\Kk}{{\cal{K}}}
\newcommand{\aaa}{{\mathbf{a}}}
\newcommand{\vv}{{\mathbf{v}}}

\newcommand{\semi}{\rtimes}

\newcommand{\nobf}{\noindent\bf}

\def\qed{\unskip\nobreak\hfil\penalty50\hskip1.75em\null\nobreak\hfil
$\blacksquare$ {\parfillskip=0pt \finalhyphendemerits=0 \par}\goodbreak}

\begin{document}
\title{The topology of the monodromy map of \\ the second order ODE}
\author{Dan Burghelea, Nicolau C. Saldanha and Carlos Tomei}
\maketitle

\begin{abstract}
We consider the following question: given $A \in SL(2,\RR)$,
which potentials $q$ for the second order Sturm-Liouville problem
have $A$ as its Floquet multiplier?
More precisely, define the monodromy map $\mu$
taking a potential $q \in L^2([0,2\pi])$
to $\mu(q) = \tilde\Phi(2\pi)$,
the lift to the universal cover $G = \widetilde{SL(2,\RR)}$ of $SL(2,\RR)$
of the fundamental matrix map $\Phi: [0,2\pi] \to SL(2,\RR)$,
\[ \Phi(0) = I, \quad \Phi'(t) = 
\begin{pmatrix} 0 & 1 \\ q(t) & 0 \end{pmatrix} \Phi(t). \]
Let $\HH$ be the real infinite dimensional separable Hilbert space:
we present an explicit diffeomorphism $\Psi: G_0 \times \HH \to H^0([0,2\pi])$
such that the composition $\mu \circ \Psi$ is the projection
on the first coordinate.
The key ingredient is the correspondence between potentials $q$
and the image in the plane of the first row of $\Phi$,
parametrized by polar coordinates, which we call the Kepler transform.
As an application among others, let $\Cc_1 \subset L^2([0,2\pi])$ be the set
of potentials $q$ for which the equation $-u'' + qu = 0$
admits a nonzero periodic solution: $\Cc_1$ is diffeomorphic to
the disjoint union of a hyperplane and 
cartesian products of the usual cone in $\RR^3$ with $\HH$.
\end{abstract}

\medbreak

{\noindent\bf Keywords:} Sturm-Liouville, monodromy, Floquet matrix,
Kepler transform.

\smallbreak

{\noindent\bf MSC-class:} 34B05; 34B24; 46T05.

\section{Introduction}
\label{section:intro}

For a given potential $q \in H^0([0,2\pi]) = L^2([0,2\pi])$,
the homogeneous equation
\begin{equation*}
- v''(t) + q(t) v(t) = 0, \quad t \in [0,2\pi]
\tag*{$(\ast)$} \label{eq:ast}
\end{equation*}
admits fundamental solutions $v_1, v_2 \in H^2([0,2\pi])$,
\[ v_1(0) = 1, \; v_1'(0) = 0, \; v_2(0) = 0, \; v_2'(0) = 1. \]
The {\it fundamental matrix}
$\Phi: [0,2\pi] \to SL(2,\RR)$ is
\[ \Phi(t) = \begin{pmatrix}
v_1(t) & v_2(t) \\ v_1'(t) & v_2'(t) \end{pmatrix} \]
and evaluation at $t = 2\pi$ obtains
the {\it Floquet multiplier} $\Phi(2\pi) \in SL(2,\RR)$.
We study the geometry of the set of potentials $q$
with given Floquet multiplier: it turns out that this set
has countably many connected components and in order to
describe them it is useful to consider the lifted version of these objects
to a covering map of $SL(2,\RR)$.

Denote by $\Pi: G = \widetilde{SL(2,\RR)} \to SL(2,\RR)$
the universal cover of the group $SL(2,\RR)$.
The {\it lifted fundamental matrix} is the continuous function
$\tilde\Phi: [0,2\pi] \to G$, $\tilde\Phi(0) = I$,
$\Pi \circ \tilde\Phi = \Phi$ and
the {\it monodromy map}, the lifted version of the Floquet multiplier,
is $\mu: H^0([0,2\pi]) \to G$, $\mu(q) = \tilde\Phi(2\pi)$.
As we shall see,
the image of $\mu$ is an open set $G_0 \subset G$ diffeomorphic to $\RR^3$.
The map $\mu$ is topologically rather simple.
Let $\HH$ be the real infinite dimensional separable Hilbert space.

\begin{theo}
\label{theo:geomonodromyintro}
There exists a diffeomorphism $\Psi: G_0 \times \HH \to H^0([0,2\pi])$
such that the composition $\mu \circ \Psi$ is the projection
on the first coordinate.
\end{theo}

Thus, the set of potentials $q$ with given monodromy $g \in G_0$
is parametrized by $\Psi(g,h)$, $h \in \HH$, and is therefore
a (topological) subspace of codimension $3$.
This theorem will be extended to other function spaces
($H^p(\Ss^1)$ and $H^p([0,2\pi])$ for $p \ge 0$)
in theorem \ref{theo:geomonodromy}.

The map $\Psi$ will be constructed explicitly via the {\it Kepler transform}.
Given a potential $q$, set $\vv: [0,1] \to \RR^2 - \{0\}$,
$\vv = (v_1, v_2)$, and let $\theta: [0,1] \to \RR$
be the continuously defined argument of $\vv$
starting with $\theta(0) = 0$, i.e., 
$\vv(t)/{|\vv(t)|} = (\cos\theta(t), \sin\theta(t))$.
It turns out that the function $\theta$ is strictly increasing
and we may therefore write
\[ \vv(t) = \sqrt{\rho(\theta(t))} \; (\cos\theta(t), \sin\theta(t)),
\quad \rho: [0,\theta_M] \to (0,+\infty), \quad \theta_M = \theta(2\pi). \]
Up to differentiability class (to be detailed in section 4),
these constructions define bijections between the following three sets:
\begin{enumerate}[(a)]
\item{$\Pp$, the set of potentials $q$;}
\item{the set $\Ff$ of {\it fundamental curves} 
$\vv: [0,2\pi] \to \RR^2 - \{0\}$ for which
$\vv(0) = (1,0)$, $\vv'(0) = (0,1)$ and $\vv(t) \wedge \vv'(t) = 1$;}
\item{the set $\Kk$ of {\it orbits}:
pairs $(\theta_M, \rho)$ where $\theta_M > 0$,
$\rho: [0,\theta_M] \to (0,+\infty)$,
$\rho(0) = 1$, $\rho'(0) = 0$ and
$\int_0^{\theta_M} \rho(\theta) d\theta = 2\pi$.}
\end{enumerate}
Luckily, monodromy is easy to handle in $\Kk$:
two potentials have the same monodromy if and only if their
orbits have the same $\theta_M$, $\rho(\theta_M)$ and $\rho'(\theta_M)$.
The level sets of $\mu$ are thus parametrized by the set of positive $\rho$'s
with prescribed behavior at endpoints and integral equal to $2\pi$.

We then proceed to apply theorem \ref{theo:geomonodromy}
to the theory of periodic Sturm-Liouville operators.
Let $\Cc \subset H^0([0,2\pi])$ be the set of potentials $q$
for which equation \ref{eq:ast} admits a periodic nontrivial solution $v$.
It is easy to see that $q \in \Cc$ if and only if $\tr\mu(q) = 2$,
thus reducing the study of $\Cc$ to the study of
the set of matrices in $G_0$ with trace equal to $2$.
The upshot is the following:
let $\Sigma_0 \subset \RR^3$ be the plane $z = 0$ and,
for $n > 0$, let $\Sigma_n$ be the cone
\[ x^2 + y^2 = \tan^2 z, \quad
2\pi n - \frac{\pi}{2} < z < 2\pi n + \frac{\pi}{2} \]
and $\Sigma = \bigcup_{n \ge 0} \Sigma_n$.

\begin{theo}
\label{theo:periodic}
There is a diffeomorphism between $(\RR^3, \Sigma) \times \HH$
and $(H^0([0,2\pi]), \Cc)$.
\end{theo}

The images of the vertices of the cones in $\Sigma \times \HH$
form a countable union of topological subspaces of codimension $3$,
the set of potentials $q$ for which all solutions of equation \ref{eq:ast}
are periodic.

Standard oscillation theory is incorporated in the following geometric
property, stated in theorem \ref{theo:oscillation}.
Consider a straight line in $H^0([0,2\pi])$
of the form $q_0 + s q_+$, $s \in \RR$,
where $q_+$ is almost everywhere strictly positive.
This line meets the image of $\Sigma_0 \times \HH$ exactly once
and the intersection is transversal.
Also, for each $n > 0$, the line meets the image of $\Sigma_n \times \HH$
either exactly twice (transversally, once in each leaf)
or once at the image of a vertex.

As an application, we describe the critical set of the nonlinear periodic
Sturm-Liouville operator with quadratic nonlinearity.
Let $p \ge 2$ and $F: H^p(\Ss^1) \to H^{p-2}(\Ss^1)$ be given by
$F(u) = - u'' + u^2/2$. 
Let $C \subset H^p(\Ss^1)$ be the critical set of $F$.
Then the pair $(H^p(\Ss^1), C)$ is diffeomorphic to
$(\RR^3, \Sigma) \times \HH$ (see corollary \ref{coro:bullet}).
This result should be contrasted to those obtained in \cite{Ruf} and
\cite{BT} for a nonlinear Sturm-Liouville operator
with Dirichlet boundary conditions and convex nonlinearity.
In \cite{BST}, the authors characterized the critical set 
with the weaker, generic hypothesis on the nonlinearity:
the components of the critical set are topological hyperplanes.
Analogous results for the periodic case,
the original motivation for this paper,
will be discussed in a forthcoming paper (\cite{BST2}).

The counterpart to the set of vertices of $\Cc$ in the third order case
is the set $C^\ast_{3,p} \subset (H^3(\Ss^1))^2$
of pairs of potentials $(q_0, q_1)$
for which {\it all} solutions $v$ of 
\[ v'''(t) - q_1(t) v'(t) - q_0(t) v(t) = 0 \]
are periodic.
Using monodromy arguments (\cite{ST}), this set is shown to be homeomorphic
to the set of closed locally convex curves in $\Ss^2$
with a prescribed basepoint,
a very complicated space with nontrivial homology for every even dimension
(\cite{Saldanha}).

The problem of characterizing potentials having $0$ in the spectrum
is clearly related to the description of isospectral classes of potentials,
as accomplished in \cite{Trubowitz}, \cite{PT} and \cite{KP}.
However, we do not think our results are
corollaries of these powerful techniques.

Back to the linear Sturm-Liouville problem, we proceed to consider
more general boundary conditions.
For a $2 \times 4$ real matrix $U$, we say a solution $v$
of equation \ref{eq:ast} satisfies $U$-boundary conditions if
\[ U \begin{pmatrix} v(0) & v'(0) & v(2\pi) & v'(2\pi) \end{pmatrix}^\ast
= \begin{pmatrix} 0 & 0 \end{pmatrix}^\ast. \]
We are again interested in the geometry and topology of $\Cc$,
the set of potentials $q$ for which equation \ref{eq:ast}
admits a nontrivial solution satisfying $U$-boundary conditions.
This again can be reduced to the study of certain algebraically defined
subsets of $G_0$.

In section 2 we present the relevant geometric facts about
$G = \widetilde{SL(2,\RR)}$ and in section 3 we do the same for
$SL^\pm(2,\RR)$, the group of real $2 \times 2$ matrices with
determinant $\pm 1$. In section 4 we present the monodromy map $\mu$
and the Kepler transform which is then used in section 5 to prove
theorem \ref{theo:geomonodromy},
a more general version of theorem \ref{theo:geomonodromyintro} above.
In section 6 we study the periodic Sturm-Liouville problem,
proving theorems \ref{theo:triplesp} and \ref{theo:oscillation},
improved versions of theorem \ref{theo:periodic}.
Finally, in section 7, we study more the Sturm-Liouville problem
with more general boundary conditions.

The last two authors received the support of CNPq, CAPES and FAPERJ (Brazil).
The second author acknowledges the hospitality of
The Mathematics Department of The Ohio State University
during the winter quarter of 2004.


\section{Coordinates for the universal cover of $SL(2,\RR)$}

Consider the universal cover
and $\Pi: G  = \widetilde{SL(2,\RR)} \to SL(2,\RR)$:
several systems of coordinates for the Lie group $G$ will be useful.
We begin with the diffeomorphism induced by the {\it Cartan decomposition}:
$\phi_C: \RR^3 \to G$ with $\phi_C(0,0,0) = I \in G$ and
\[ (\Pi \circ \phi_C)(\alpha, r \cos \eta, r \sin \eta) = \]
\begin{equation} \begin{pmatrix} \cos\alpha & \sin\alpha \\
-\sin\alpha & \cos\alpha \end{pmatrix}
\begin{pmatrix} \cosh r + \sinh r \cos\eta & \sinh r \sin\eta \\
\sinh r \sin\eta & \cosh r - \sinh r \cos\eta \end{pmatrix}. 
\label{eq:Cartan} \end{equation}


We are interested in the stratification of $G$ in conjugacy classes.
The center $Z(G)$ of $G$ is formed by the elements of the form $\iota^n$,
$n \in \ZZ$, where $\iota = \phi_C(\pi,0,0)$:
we have $\Pi(\iota^n) = (-1)^n I$.
From the connectivity of $G$, conjugacy classes are contained in
connected components of level sets $T_c = \tr^{-1}(\{c\})$
of the trace function $\tr: G \to \RR$.
We systematically abuse notation by writing $\tr g$ instead of $\tr(\Pi g)$.
For any matrix $A \in SL(2,\RR)$,
$A \ne \pm I$, the centralizer $\{ B \;|\; AB = BA \}$ is a
Lie group of dimension $1$ and, since $G$ is a covering of $SL(2,\RR)$,
the same holds for the centralizer of any $g \in G$,
$g \ne \iota^n$, $n \in \ZZ$.
Thus, the conjugacy class of any such $g$ is a $2$-dimensional manifold.

\begin{figure}[ht]
\begin{center}
\psfrag{pi/2}{$\pi/2$}
\psfrag{pi}{$\pi$}
\psfrag{3pi/2}{$3\pi/2$}
\psfrag{-pi/2}{$-\pi/2$}
\psfrag{-pi}{$-\pi$}
\psfrag{-3pi/2}{$-3\pi/2$}
\psfrag{alpha}{$\alpha$}
\psfrag{r}{$r$}
\psfrag{A-1}{$A_{-1}$}
\psfrag{A0}{$A_0$}
\psfrag{A1}{$A_1$}
\psfrag{A2}{$A_2$}
\epsfig{height=44mm,file=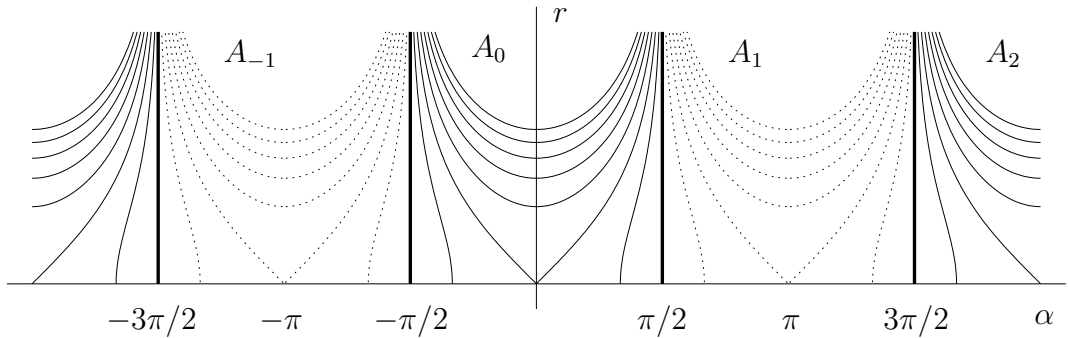}
\end{center}
\caption{Level curves of the trace function}
\label{fig:trlevel}
\end{figure}

A straightforward computation yields
$\tr\phi_C(\alpha, r \cos \eta, r \sin \eta) = 2 \cos \alpha \cosh r$. 
The sets $\phi_C^{-1}(T_c)$ are obtained by rotating
figure \ref{fig:trlevel} around the horizontal axis ($r = 0$).
The figure indicates the level curves for $c \in \ZZ$,
solid for $c > 0$, thicker for $c = 0$ and dotted for $c < 0$.
The V shaped curves correspond to $c = \pm 2$.
Notice that $T_0$ is the countable union of planes
$\alpha = k\pi + \pi/2$ in Cartan coordinates.

The sign of the trace is determined by $\cos\alpha$.
Defining $A_n = \phi_C((n\pi - \pi/2, n\pi + \pi/2)\times\RR^2)$,
the regions bounded by the thick vertical lines in the figure,
the sign of the trace is constant equal to $(-1)^n$ in each open
set $A_n$. Since $A_n = \iota^n A_0$ it suffices to study
the trace function in $A_0$.
From the picture, level sets $T_c$ look like cones or hyperboloids.
To make this precise, define the real analytic functions
\[ f_1(x) = \frac{\arccos(\exp(-x^2))}{|x|}, \quad
f_2(x) = \frac{\arccosh(\exp(x^2))}{|x|} \]
and $\phi_X: \RR^3 \to A_0$ by
\[ \phi_X(x,y,z) =
\phi_C\left(x f_1(x), y f_2(\sqrt{y^2+z^2}), z f_2(\sqrt{y^2+z^2})\right): \]
it is easy to verify that $\phi_X$ is a diffeomorphism and that
\[ \tr(\phi_X(x,y,z)) = 2\exp(-x^2+y^2+z^2). \]
Thus, for $c > 0$, $\phi_X^{-1}(T_c \cap A_0)$ is
the surface $-x^2 + y^2 + z^2 = \log(c/2)$.
For $0 < c < 2$ this is a hyperboloid with two connected components,
diffeomorphic to the disjoint union of two planes:
in this case, the set $T_c$ is a disjoint union
of countably many surfaces diffeomorphic to $\RR^2$, two in each $A_{2n}$.
For $c > 2$, $\phi_X^{-1}(T_c \cap A_0)$ is a one-sheet hyperboloid,
diffeomorphic to $\Ss^1 \times \RR$:
here, $T_c$ is a disjoint union of countably many surfaces
diffeomorphic to $\Ss^1 \times \RR$, one in each $A_{2n}$.
Finally, $\phi_X^{-1}(T_2 \cap A_0)$ is the cone $x^2 = y^2 + z^2$,
which, except for one point, the {\it vertex}, is a submanifold.
We call the cone $\cone$.
Thus, for $c = 2$, $T_c$ is a disjoint union
of countably many copies of $\cone$, one in each $A_{2n}$
The connected component of $T_2$ containing $I$
is the image under the exponential map of the cone of nilpotent matrices
in the Lie algebra of $G$ (naturally identified with $sl(2,\RR)$).
The cases $c < -2$, $c = -2$ and $-2 < c < 0$ are similar,
with the components now lying in $A_{2n+1}$.

Summing up, for each $c \ne \pm 2$, the connected components of $T_c$
are conjugacy classes in $G$.
The vertices of the cones in $T_{\pm 2}$ are precisely $\iota^n$:
each vertex is a conjugacy class by itself.
A cone minus the vertex consists of two {\it leaves},
each of them diffeomorphic to $\Ss^1 \times \RR$:
each leaf of a cone is a conjugacy class.
Let $T_2^0 \subset T_2$ be the connected component containing the origin.
The two leaves of the cone $T_2^0$ meet at the vertex $I$
and consist of lifted matrices with both eigenvalues equal to $1$.
Thus, $g \in T_2^0 - \{I\}$ projects to $I+N \in SL(2,\RR)$,
$N$ a nonzero nilpotent matrix.
Define $\sgn(g)$ to be $\sgn(\det(Nv,v))$, $v \notin \ker N$;
this sign is well defined and may be used as a label for the leaf.

Consider now the {\it left Iwasawa decomposition}
$\phi_L: \RR \times (0,\infty) \times \RR \to G$ with
$\phi_L(0,1,0) = I$ and
\begin{equation}
(\Pi \circ \phi_L)(\theta,\rho,\nu) =
\begin{pmatrix} \cos\theta & \sin\theta \\
-\sin\theta & \cos\theta \end{pmatrix}
\begin{pmatrix} \sqrt{\rho} & 0 \\ 0 & 1/\sqrt{\rho} \end{pmatrix} 
\begin{pmatrix} 1 & 0 \\ \nu/2 & 1 \end{pmatrix}
\label{eq:leftIwasawa} \end{equation}
and the open nested half-spaces
$G_{\theta} = \phi_L((\theta,+\infty) \times (0,+\infty) \times \RR)
\subset G$.
The set $G_{\theta}$ consists of the elements $g \in G$
for which the variation in argument from $e_2$ to $g e_2$
is smaller than $-\theta$
(the variation in argument is computed along a path
$\gamma: [0,1] \to G$ joining $\gamma(0) = I$ to $\gamma(1) = g$).
The pairs $(G_\theta, T_c \cap G_\theta)$ come up in the study
of the monodromy map in later sections.

\begin{prop}
\label{prop:fivepairs}
For any $\theta$ and $c$, the pair $(G_\theta, T_c \cap G_\theta)$ 
is diffeomorphic to $(G_0, T_{\hat c} \cap G_0)$
for $\hat c = 0$, $\hat c = \pm 2$ or $\hat c = \pm 4$.
More precisely,
\[ \hat c = \begin{cases} 0,&|c| < 2,\\
(-1)^{\lfloor \theta/\pi \rfloor} c,&|c| = 2, \\
4 (-1)^{\lfloor \theta/\pi \rfloor} \sgn(c),&|c| > 2. \end{cases} \]
\end{prop}

Recall that $\lfloor x \rfloor$ is the only integer in the interval $(x-1,x]$.
Along the proof of the proposition, we will give geometric descriptions
of the five pairs in the statement.

\begin{figure}[ht]
\begin{center}
\psfrag{T20}{$T_2^0$}
\psfrag{T0}{$T_0$}
\psfrag{T2}{$T_2$}
\psfrag{T-2}{$T_{-2}$}
\psfrag{T4}{$T_4$}
\psfrag{T-4}{$T_{-4}$}
\psfrag{partialG0}{$\partial G_0$}
\epsfig{height=30mm,file=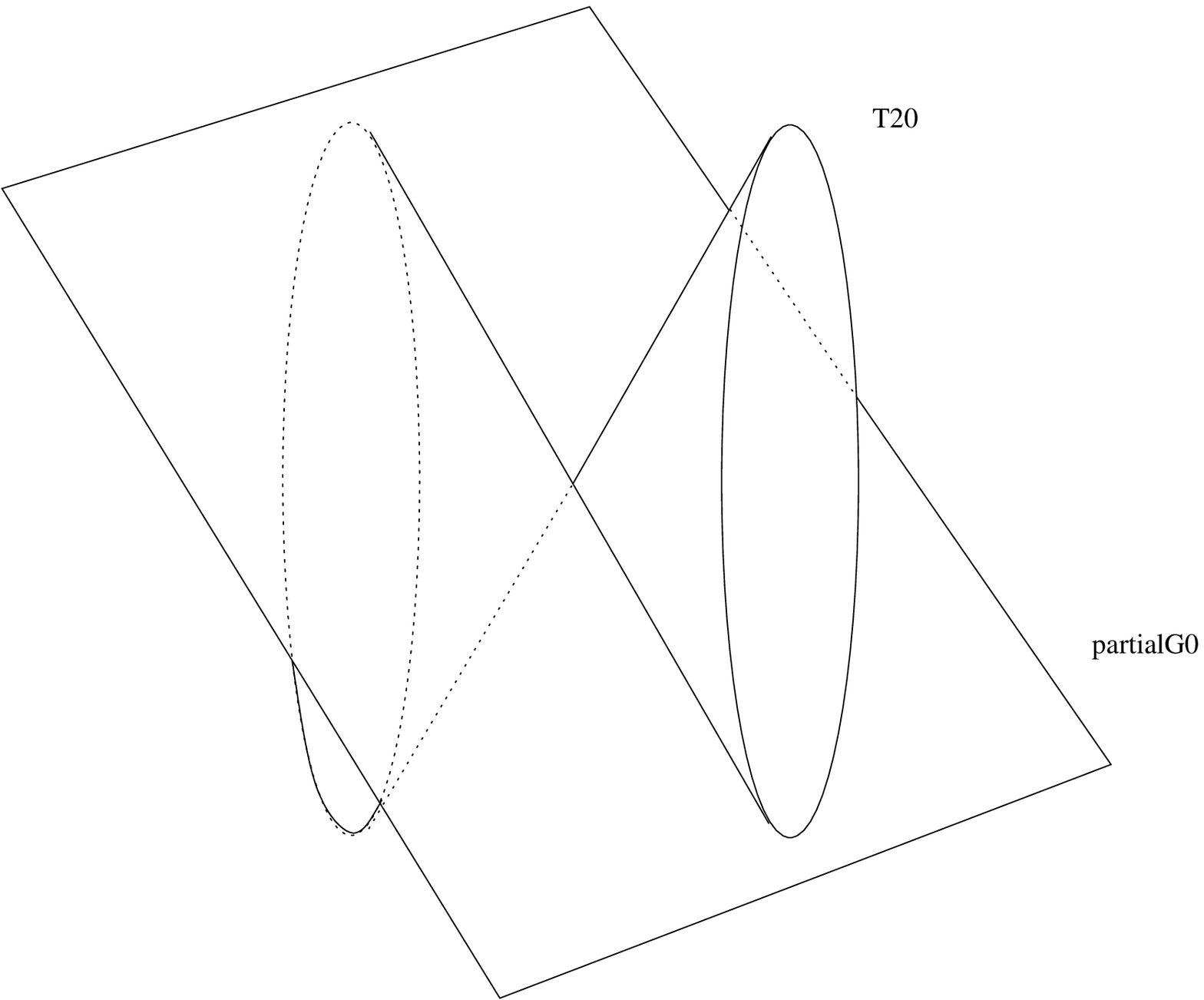}
\qquad
\epsfig{height=40mm,file=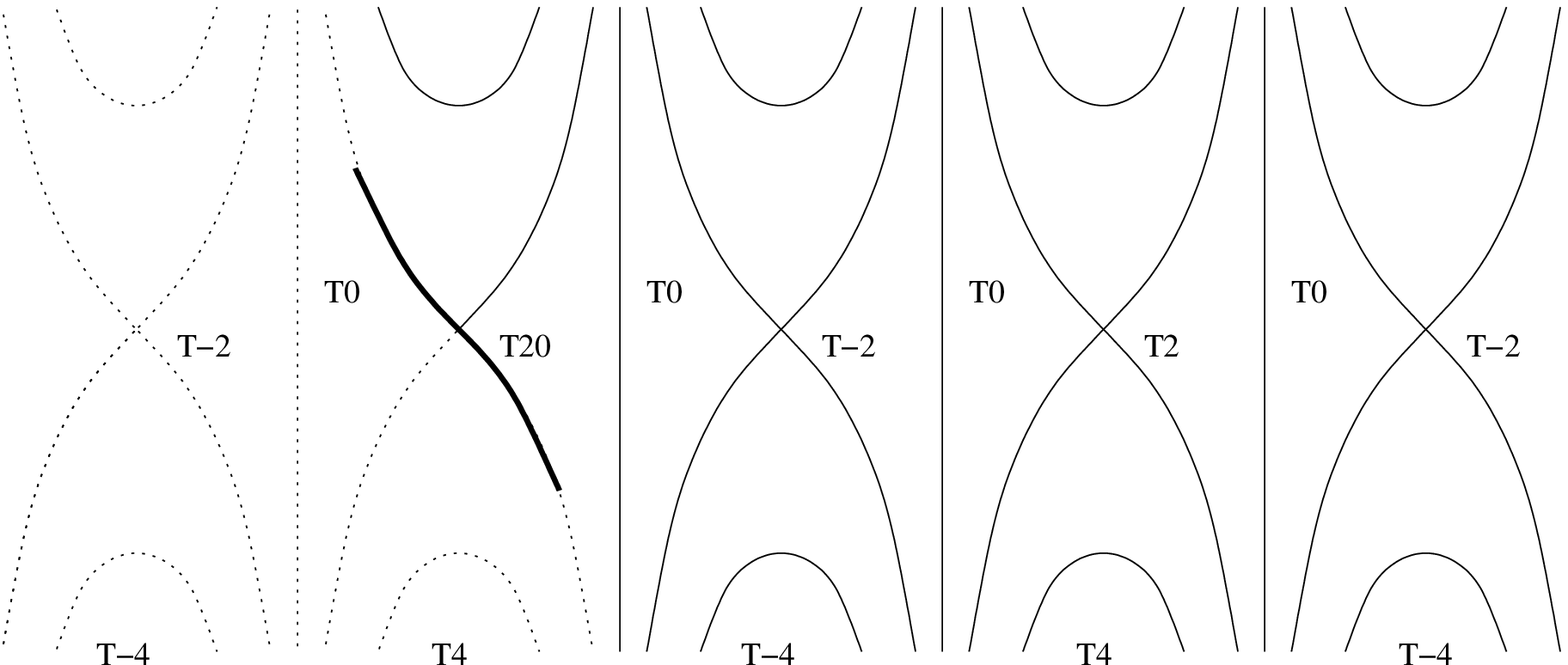}
\end{center}
\caption{Two views of $G_0 \subset G$.}
\label{fig:g0}
\end{figure}

{\nobf Proof:}
Since $\phi_L$ is a diffeomorphism,
the boundaries $\partial G_\theta$ are smooth (topological) hyperplanes.
The surface $\partial G_0$ consists of (lifts of) lower triangular matrices
with positive diagonal entries.
Clearly, for $g \in \partial G_0$, $\tr g \ge 2$,
and on the curve of lower triangular matrices with diagonal $(1,1)$
we have $\tr g = 2$.
This implies that the surface $\partial G_0$ is tangent to the cone $T_2^0$.
For $g \in T_2^0$, except for the curve of tangency,
$\sgn(g)$ coincides with the sign of $\theta$:
indeed, the sign $\sgn(g)$ is also the sign of the variation of argument
from $gv$ to $v$ if $v$ is not an eigenvector of $g$.
Thus, the positive leaf of $T_2^0$ is contained
in the closure of $G_0$ and the negative one is disjoint from $G_0$.
The intersection of $T_2^0$ with $G_0$ is therefore
the positive leaf minus a closed half-line: 
it is thus diffeomorphic to a plane.
Figure \ref{fig:g0} shows the set $G_0$, together with the cones $T_{\pm 2}$,
in two kinds of representations. The drawing on the left is an attempt
to give a 3d perspective view of $T_2^0$ and $\partial G_0$
as a cone and a tangent plane.
The drawing on the right is far more schematic:
the connected components of $T_2$ and $T_{-2}$ are shown as big Xs,
the parts contained in $G_0$ drawn in solid lines
and the others in dotted lines;
$\partial G_0$ is represented by a thick line.

The connected component of $T_2 \cap G_0$ contained in $A_0$,
the solid half-line starting at the thick line in figure \ref{fig:g0},
is (diffeomorphic to) a plane while the other components,
one in each $A_{2k}$, $k > 0$,
are (diffeomorphic to) cones with horizontal axis.
On the other hand, the components of $T_{-2} \cap G_0$,
one in each $A_{2k+1}$, $k \ge 0$, are all cones.
In particular, the pairs $(G_0, T_2 \cap G_0)$ and $(G_0, T_{-2} \cap G_0)$
are not diffeomorphic.

Similarly, as we shall soon prove, the connected component of $T_4 \cap G_0$
in $A_0$, drawn as a branch of a fake hyperbola in figure \ref{fig:g0},
is a plane, while the other components, drawn as complete fake hyperbola,
are cylinders (i.e., diffeomorphic to $\Ss^1 \times \RR$).
The components of $T_{-4} \cap G_0$ are cylinders
and those of $T_0 \cap G_0$, drawn as vertical lines, are planes.

Let $B_n \subset G$ be $G_{n\pi} - \bar G_{(n+1)\pi} =
\phi_L((n\pi, (n+1)\pi) \times (0,+\infty) \times \RR)$.
Thus, the sets $B_n$ are open and disjoint and, together with the sets $A_n$,
form an open cover of $G$ with $A_n \cap B_{n'} \ne \emptyset$
if and only if $n = n'$ or $n = n'+1$.
The map $\phi_X$ provided a normal form for the trace on $A_n$.
For $B_n$ instead, consider the diffeomorphism
$\phi_Y: (n\pi,(n+1)\pi) \times (0,+\infty) \times \RR \to B_n \subset G$, 
\[ \phi_Y(\theta,\rho,c) = \phi_L\left(\theta,\rho,
\frac{2\rho^{1/2}(c-(\rho^{1/2}+\rho^{-1/2})\cos\theta)}{\sin\theta}
\right), \]
for which $\tr(\phi_Y(\theta,\rho,c)) = c$.
Thus, $(B_n, T_c \cap B_n)$ is diffeomorphic to the pair $(\RR^3, \{z=c\})$
and so is $(T_c \cap B_n \cap G_\theta, B_n \cap G_\theta)$
assuming $n\pi < \theta < (n+1)\pi$.

Consider now arbitrary values of $\theta$ and $c$.
We may assume $\theta \in [0,\pi)$
by multiplying everything in sight by an appropriate element
$\iota^{n}$ of the center of $G$,
an operation which, up to sign, preserves traces.
Set $\epsilon > 0$, $\epsilon < \pi - \theta$.
The diffeomorphism $\phi_Y$ yields a diffeomorphism between
the regions $G_0 - G_{\theta+\epsilon}$ and $G_\theta - G_{\theta+\epsilon}$,
coinciding with the identity near their common boundary and preserving trace.
We therefore have a diffeomorphism between the pairs
$(G_0, T_c \cap G_0)$ and $(G_\theta, T_c \cap G_\theta)$,
which, together with the geometric descriptions in figure \ref{fig:g0},
completes the proof.
\qed


\section{$SL^-(2,\RR)$}

Let $SL^{\pm}(2,\RR) \subset GL(2,\RR)$ be the group of matrices
of determinant $\pm 1$. Clearly, $SL^{\pm}(2,\RR)$ has two connected
components: $SL^+(2,\RR) = SL(2,\RR)$ and $SL^-(2,\RR)$,
the set of $2 \times 2$ real matrices of determinant $-1$. For 
\[ H = \{ I, R \}, \quad
R = \begin{pmatrix} -1 & 0 \\ 0 & 1 \end{pmatrix}, \]
$SL^\pm(2,\RR)$ is the semidirect product $SL(2,\RR) \semi H$
defined by the outer automorphism $r: SL(2,\RR) \to SL(2,\RR)$,
\[ r\left( \begin{pmatrix} a & b \\ c & d \end{pmatrix}\right) =
R \begin{pmatrix} a & b \\ c & d \end{pmatrix} R =
\begin{pmatrix} a & -b \\ -c & d \end{pmatrix} . \]
The automorphism $r$ lifts to $\tilde r: G \to G$,
also an automorphism of order $2$.
Use $\tilde r$ to define $G^\pm$ as the semidirect product $G \semi (\ZZ/(2))$.
More concretely, set $G^\pm$ to be the disjoint union of $G$ and
$\tilde R G = \{ \tilde R g, g \in G\}$, the product being defined by
$g\tilde R = \tilde R \tilde r(g)$ and $\Pi^\pm: G^\pm \to SL^\pm(2,\RR)$
is a homomorphism extending $\Pi: G \to SL(2,\RR)$ with $\Pi(\tilde R) = R$.
Clearly, $G^\pm$ has two connected components $G^+ = G$ and $G^-$,
each homeomorphic to $\RR^3$ and the projection $\Pi^\pm$ is
a universal cover on each connected component.

The {\it Schur decomposition} induces a diffeomorphism
$\phi_S: \RR \times (0,+\infty) \times \RR \to G^-$
with $\phi_S(0,1,0) = \tilde R$,
\begin{equation}
(\Pi \circ \phi_S)(\alpha, \lambda, \nu) =
\begin{pmatrix} \cos\alpha & \sin\alpha \\
-\sin\alpha & \cos\alpha \end{pmatrix}
\begin{pmatrix} -1/\lambda & 0 \\ \nu & \lambda \end{pmatrix}
\begin{pmatrix} \cos\alpha & -\sin\alpha \\
\sin\alpha & \cos\alpha \end{pmatrix}.
\label{eq:Schur} \end{equation}
As before, we consider the level sets $T^-_c = \tr^{-1}(c) \subset G^-$.
Clearly, $T^-_c = \phi_S(\RR \times \{\lambda\} \times \RR)$
where $\lambda$ is the (only) positive solution of $\lambda - 1/\lambda = c$.
This implies that $T^-_c$ is always (diffeomorphic to) a plane.

The {\it left Iwasawa decomposition} is the diffeomorphism
$\phi^-_L: \RR \times (0,\infty) \times \RR \to G^-$,
$\phi^-_L(\theta,\rho,\nu) = \tilde R \phi_L(\theta,\rho,\nu)$.
Finally, set
$G^-_\theta = \phi^-_L((\theta,+\infty)\times (0,\infty) \times \RR)$.
The topology of the pairs $(G^-_\theta, T^-_c \cap G^-_\theta)$
is much simpler than that of their positive counterparts.

\begin{prop}
\label{prop:onepair}
For any $\theta$ and $c$, the pair $(G^-_\theta, T^-_c \cap G^-_\theta)$
is diffeomorphic to $(\RR^3, \{z = 0\})$.
\end{prop}

{\nobf Proof:}
For $\lambda > 0$ and $\nu \in \RR$,
define $\theta_{\lambda,\nu}: \RR \to \RR$,
so that $\theta_{\lambda,\nu}(\alpha)$ is the first coordinate of
$(\phi^-_L)^{-1}(\phi_S(\alpha,\lambda,\nu))$.
A straightforward computation verifies that
$\theta_{\lambda,\nu}$ is a diffeomorphism:
informally, given $\lambda$ and $\nu$,
the variables $\alpha$ and $\theta$ are interchangeable.
In other words, there exists a diffeomorphism
$\phi_Z: \RR \times (0,+\infty) \times \RR \to G^-$
such that $(\phi^-_L)^{-1}(\phi_Z(\theta,\lambda,\nu)) = (\theta, \ast, \ast)$
and $(\phi_S)^{-1}(\phi_Z(\theta,\lambda,\nu)) = (\ast, \lambda, \nu)$.
The result is now obvious.
\qed

\section{Monodromy and the Kepler transform}

Let $H^p([0,2\pi])$, $p \ge 0$, be the Sobolev space
of real functions whose $p$-th derivative is
in $H^0([0,2\pi]) = L^2([0,2\pi])$.
We also consider the {\it periodic} Sobolev space
$H^p(\Ss^1) \subset H^p([0,2\pi])$ of $u$'s
with $u(0) = u(2\pi), \ldots, u^{(p-1)}(0) = u^{(p-1)}(2\pi)$.
The periodic space $H^p(\Ss^1)$ is a closed subspace
of $H^p([0,2\pi])$ of codimension $p$.

For a given potential $q \in H^p([0,2\pi])$, $p \ge 0$,
the {\it fundamental solutions} of the homogeneous equation
\begin{equation*}
- v''(t) + q(t) v(t) = 0, \quad t \in [0,2\pi]
\tag*{$(\ast)$} 
\end{equation*}
are those functions $v_i \in H^{p+2}([0,2\pi])$ with initial conditions
\[ v_1(0) = 1, v_1'(0) = 0, \qquad v_2(0) = 0, v_2'(0) = 1. \]
Equivalently, write 
\[ \Phi(t) = \begin{pmatrix} v_1(t) & v_2(t) \\
v_1'(t) & v_2'(t) \end{pmatrix}\]
so that
\[ \Phi(0) = I, \quad \Phi'(t) = 
\begin{pmatrix} 0 & 1 \\ q(t) & 0 \end{pmatrix} \Phi(t). \]
The fact that the Wronskian of $v_1$ and $v_2$ is constant equal to $1$
implies that $\Phi(t) \in SL(2,\RR)$ for all $t$:
thus $\Phi$ is a continuous (actually, $H^{p+1}$) function
from $[0,2\pi]$ to $SL(2,\RR)$.
Define the lifted fundamental matrix $\tilde\Phi: [0,2\pi] \to G$
by $\tilde\Phi(0) = I$ and $\Phi = \Pi \circ \tilde\Phi$
where $\Pi$ is the natural projection from $G$ to $SL(2,\RR)$.
Any solution of the homogeneous equation \ref{eq:ast} is of the form
\[ v(t) = \begin{pmatrix} 1 & 0 \end{pmatrix} 
\Phi(t) \begin{pmatrix} a_1 \\ a_2 \end{pmatrix} \]
for real constants $a_1$ and $a_2$.
Define the {\it monodromy} $\mu: H^0([0,2\pi]) \to G$,
$\mu(q) = \tilde\Phi(2\pi)$;
thus, $\mu(q)$ contains (discretely) more information than
the Floquet multiplier $\Phi(2\pi) = \Pi(\mu(q))$.

We now construct smooth natural bijections between the following three sets:
\begin{enumerate}[(a)]
\item{$\Pp = H^p([0,2\pi])$, the set of potentials $q$;}
\item{the set $\Ff$ of {\it fundamental curves}:
paths $\vv: [0,2\pi] \to \RR^2 - \{0\}$ of class $H^{p+2}$ satisfying
$\vv(0) = (1,0)$, $\vv'(0) = (0,1)$ and $\vv(t) \wedge \vv'(t) = 1$
for all $t$;}
\item{the set $\Kk$ of {\it orbits}:
pairs $(\theta_M, \rho)$ where $\theta_M > 0$
is a real number and $\rho: [0,\theta_M] \to (0,+\infty)$
is a function of class $H^{p+2}$ satisfying
$\rho(0) = 1$, $\rho'(0) = 0$ and
$\int_0^{\theta_M} \rho(\theta) d\theta = 2\pi$.}
\end{enumerate}

Let $\vv: [0, 2\pi] \to \RR^2 - \{0\}$ be the first row of $\Phi$:
thus, $\vv$ is a continuous function satisfying
\[ \vv''(t) = q(t) \vv(t), \quad \vv(0) = (1,0), \quad \vv'(0) = (0,1). \]
The condition $\det\Phi = 1$ is translated as $\vv \wedge \vv' = 1$:
in particular, the argument $\theta$ of $\vv$ always has positive derivative.
We call $\vv$ the fundamental curve associated with the potential $q$:
the map from $\Pp$ to $\Ff$ takes $q$ to $\vv$.

This map is indeed a continuous bijection:
if $\vv \in \Ff$, $\vv$ of class $H^{p+2}$,
we have $\vv \wedge \vv' = 1$ so that $\vv \wedge \vv'' = 0$.
Since $\vv$ is continuous and nonzero, $\vv''$ is a multiple of $\vv$,
i.e., $\vv'' = q \vv$ and it is straightforward to check that 
\begin{equation}
q(t) = \vv''(t) \wedge \vv'(t) \label{eq:hfromvv}
\end{equation}
and the potential $q$ lies in $H^{p}([0,2\pi])$
with $\vv$ being its associated fundamental curve.

Let $\theta: [0,2\pi] \to \RR$ be the continuously defined
argument of $\vv = (v_1, v_2)$ with $\theta(0) = 0$.
The condition $\vv \wedge \vv' = 1$ indicates that the area surrounded
by the curve $\vv$ in an interval $[t_1, t_2]$ is $(t_2 - t_1)/2$
and therefore the argument $\theta$ is strictly increasing.
Set $\theta_M = \theta(2\pi)$ and consider
$\rho: [0,\theta_M] \to (0,+\infty)$ and $\nu: [0, \theta_M] \to \RR$
defined by $\rho(\theta(t)) = |\vv(t)|^2$ and 
$\nu(\theta) = \rho'(\theta)/\rho(\theta)$.
Notice that
$t_2 - t_1 =  \int_{\theta(t_1)}^{\theta(t_2)} \rho(\theta) d\theta$
whence, in particular, $\int_0^{\theta_M} \rho(\theta) d\theta = 2\pi$.
We just constructed the map from $\Ff$ to $\Kk$;
the conditions $\rho(0) = 1$, $\rho'(0) = 0$ are easy to check.

The values of $\theta$, $\rho$ and $\nu$ admit an interpretation
in terms of the {\it right Iwasawa decomposition}.
Define the diffeomorphism
$\phi_R: \RR \times (0,\infty) \times \RR \to G$
by $\phi_R(0,1,0) = I$ and
\begin{equation}
(\Pi \circ \phi_R)(\theta,\rho,\nu) =
\begin{pmatrix} \sqrt{\rho} & 0 \\ 0 & 1/\sqrt{\rho} \end{pmatrix} 
\begin{pmatrix} 1 & 0 \\ \nu/2 & 1 \end{pmatrix}
\begin{pmatrix} \cos\theta & \sin\theta \\
-\sin\theta & \cos\theta \end{pmatrix}.
\label{eq:rightIwasawa} \end{equation}
It is easy to verify that
$\tilde\Phi(t) = \phi_R(\theta(t), \rho(\theta(t)), \nu(\theta(t)))$.

\begin{figure}[ht]
\begin{center}
\epsfig{height=35mm,file=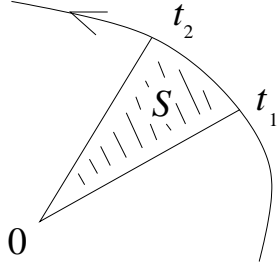}
\end{center}
\caption{$S = (t_2 - t_1)/2$:
the curve $\vv$ sweeps equal areas in equal times.}
\label{fig:kepler}
\end{figure}

The orbit $(\theta_M, \rho)$ yields a curve in the plane.
We can uniquely parametrize it so that it sweeps area $t/2$ in time $t$,
turning the curve into a fundamental curve and
thus constructing the inverse map from $\Kk$ to $\Ff$:
\begin{align}
\theta'(t) &= \frac{1}{\rho(\theta(t))},
\label{eq:thetaprimefromrho}
\\
q(t) &= \left(\frac{2 \rho'' \rho - 3 (\rho')^2 - 4 \rho^2}{4 \rho^4}\right)
(\theta(t)),
\label{eq:qfromrho}
\\
q'(t) &= \left(\frac{\rho''' \rho^2 - 7 \rho'' \rho' \rho + 6(\rho')^3
+ 4\rho' \rho^2}{2 \rho^6}\right)(\theta(t)).
\label{eq:qprimefromrho}
\end{align}
The fact that these bijections preserve smoothness class
is left to the reader.
We call this bijection between $\Ff$ and $\Kk$ the {\it Kepler transform}.

The restrictions of these bijections to the periodic case work well
but, for $p > 0$, we still have to describe the image in $\Ff$ and $\Kk$
of $H^{p}(\Ss^1) \subset H^{p}([0,2\pi]) = \Pp$.
More precisely, we translate the conditions
$q^{(j)}(0) = q^{(j)}(2\pi)$, $0 \le j < p$,
in terms of the functions $\vv$ and $\rho$.
For $\vv$, we clearly must have
$\vv^{(j)}(2\pi) = \vv^{(j)}(0) \mu(q)$, $2 \le j < p+2$.
For $\rho$, the conditions become far more complicated.
From equations \ref{eq:qfromrho} and \ref{eq:qprimefromrho},
the conditions $q(0) = q(2\pi)$ and $q'(0) = q'(2\pi)$ become
\begin{align}
\rho''(\theta_M) &=
(\rho(\theta_M))^3 \rho''(0) + b_0(\rho(\theta_M),\rho'(\theta_M)) \notag\\
\rho'''(\theta_M) &=
(\rho(\theta_M))^4 \rho'''(0) + b_1(\rho(\theta_M),\rho'(\theta_M),\rho''(0))
\notag\end{align}
where $b_0$ and $b_1$ are smooth functions.
More generally, formulae for higher derivatives of $q$
yield a translation from $q^{(j)}(0) = q^{(j)}(2\pi)$ to 
\[ \rho^{(j+2)}(\theta_M) =
(\rho(\theta_M))^{j+3} \rho^{(j+2)}(0) +
b_j(\rho(\theta_M),\rho'(\theta_M),\rho''(0),\ldots,\rho^{j+1}(0)), \]
where $b_j$ is a rather complicated expression.
Summing up, there exists smooth maps
$B_p: (0,+\infty) \times \RR \to \Diff(\RR^p,\RR^p)$
such that for $\rho \in H^{p+2}([0,\theta_M])$,
the associated potential $q \in H^{p}([0,2\pi])$
belongs to $H^{p}(\Ss^1)$ if and only if
\begin{equation}
(\rho''(\theta_M), \ldots,  \rho^{(p+1)}(\theta_M)) =
B_p(\rho(\theta_M), \rho'(\theta_M)) \; (\rho''(0), \ldots, \rho^{(p+1)}(0)).
\label{eq:beta} \end{equation}

\begin{prop}
\label{prop:imagemu}
For any $p \ge 0$, the image of $\mu: H^p([0,2\pi]) \to G$
or $\mu: H^p(\Ss^1) \to G$ is $G_0$.
\end{prop}

{\nobf Proof: }
First notice that
$\phi_L(\theta_L, \rho_L, \nu_L) = \phi_R(\theta_R, \rho_R, \nu_R)$
then $\sgn(\theta_L) = \sgn(\theta_R)$, so that
\[ G_0 = \phi_L((0,+\infty) \times (0,+\infty) \times \RR)
= \phi_R((0,+\infty) \times (0,+\infty) \times \RR). \]

Clearly, for any $q \in H^p$,
since $\theta: [0,2\pi] \to \RR$ is strictly increasing
with $\theta(0) = 0$ then $\theta_M = \theta(2\pi) > 0$
and $\mu(q) = \phi_R(\theta_M, \rho(\theta_M), \nu(\theta_M)) \in G_0$.

Conversely, take $p \in G_0$.
Write $p = \phi_R(\theta_M, \rho_M, \nu_M)$, $\theta_M > 0$.
Construct an $H^{p+2}$ function $\rho: [0, \theta_M] \to (0, +\infty)$
with $\rho(0) = 1$, $\rho'(0) = 0$, $\rho(\theta_M) = \rho_M$,
$\rho'(\theta_M) = \rho_M \nu_M$ and
\[ \int_0^{\theta_M} \rho(\theta) d\theta = 2\pi. \]
Apply the Kepler transform on the pair $(\theta_M, \rho)$
to obtain a potential $h$ with $\mu(q) = p$.
Minor adjustments at the boundary points may be performed to
guarantee that $q \in H^p(\Ss^1)$.
\qed

\section{Global geometry of the monodromy map}

We are ready to prove the first main result of this paper.
Geometrically, the theorem states that level sets of the monodromy
map are, after a smooth change of variables,
parallel affine subspaces of codimension $3$.
The claim holds for the restriction of the monodromy
to $H^p([0,2\pi])$ and to $H^p(\Ss^1)$, $p \ge 0$.
Let $\HH$ be the real separable infinite dimensional Hilbert space.

\begin{theo}
\label{theo:geomonodromy}
For $p \ge 0$, there exists
smooth diffeomorphisms
$\Psi_{[0,2\pi]}^p: G_0 \times \HH \to H^{p}([0,2\pi])$ and
$\Psi_{\Ss^1}^p: G_0 \times \HH \to H^{p}(\Ss^1)$
such that both compositions $\mu \circ \Psi^p$ are projections
on the first coordinate.
\end{theo}

The subscript $[0,2\pi]$ or $\Ss^1$ for the diffeomorphisms $\Psi^p$
will be omitted whenever it is clear from the context.
The proof yields an explicit construction of the maps $\Psi^p$.

{\nobf Proof: }
We first consider the case $p = 0$.
Take $g_0 \in G_0$ and consider its right Iwasawa coordinates
$(\theta_0, \rho_0, \nu_0) \in (0,+\infty) \times (0,+\infty) \times \RR$.
A potential $q \in H^0([0,2\pi])$ has monodromy $g_0$
if and only if its associated orbit $(\theta_M, \rho)$
(where $\rho \in H^2([0,\theta_M])$ with $\rho(\theta) > 0$
for all $\theta$, $\rho(0) = 1$, $\rho'(0) = 0$
and $\int_0^{\theta_M} \rho(\theta) d\theta = 2\pi$)
satisfies $\theta_M = \theta_0$, $\rho(\theta_M) = \rho_0$,
$\rho'(\theta_M) = \nu_0 \rho_0$.
We shall parametrize the set of all such functions $\rho$
by a Hilbert space $H = \HH$.

We first choose a base point $\Psi^0(g_0,0)$.
There exists a unique polynomial $P_0 = P_{\theta_0, \rho_0, \nu_0}$
of degree $4$ or less such that
\[ (\exp \circ P_0)(0) = 1, \quad (\exp \circ P_0)'(0) = 0, \]
\[ (\exp \circ P_0)(\theta_M) = \rho_0, \quad
(\exp \circ P_0)'(\theta_M) = \nu_0 \rho_0, \]
\[ \int_0^{\theta_M} (\exp \circ P_0)(\theta) d\theta = 2\pi. \]
The exponential is used to guarantee the positivity of the function
$\rho = \exp \circ P_0$.
Indeed, from Lagrange interpolation there exists a unique polynomial
$P_1$ of degree at most $3$ satisfying the boundary conditions;
thus, a polynomial $P$ of degree at most $4$ satisfies the
boundary conditions if and only if $P$ is of the form
$P(\theta) = P_1(\theta) + c \theta^2(\theta_M - \theta)^2$.
The integral on the fifth condition is now a continuous strictly increasing
function of $c$ ranging from $0$ to $+\infty$ as $c$ varies in $\RR$:
there exists therefore a unique value of $c$ for which $P_0 = P$
satisfies boundary and integral conditions.
Set $\Psi^0(g_0,0)$ to be the potential associated
to the orbit $(\theta_M, \exp \circ P_0)$.

Now let $H \subset H^2([0,1])$ be the closed subspace of functions $r$ with 
\[ r(0) = r'(0) = r(1) = r'(1) = \int_0^1 r(t) dt = 0. \]
Define $\Psi^0(g_0,r)$ to be the potential
with orbit $(\theta_M, \rho)$ where
\[ \rho(\theta) = \exp\left(P(\theta) + 
r(\theta/\theta_M) +
c \theta^{2}(\theta_M - \theta)^{2}\right), \]
the parameter $c$ being again uniquely chosen
so that $\rho$ satisfies the integral condition.


The nonperiodic case for $p > 0$ is similar.
We now consider the periodic case for $p > 0$.
Take $g_0 = \phi_R(\theta_M, \rho_0, \nu_0) \in G_0$.
Let $H_1 \subset H^{p+2}([0,1])$ be the space of functions $r$
for which
\[ r(0) = r(1) = r'(0) = r'(1) = \cdots = r^{(p+1)}(0) = r^{(p+1)}(1) =
\int_0^1 r(t) dt = 0 \]
and $H = \RR^{p} \times H_1$.
Let $a_0 = 1$, $a_1 = 0$, $b_0 = \rho_0$, $b_1 = \nu_0 \rho_0$.
For each $\aaa = (a_2, \ldots, a_{p+1}) \in \RR^{p}$, let
$(b_2, \ldots, b_{p+1}) = B_p(b_0,b_1)(\aaa)$
(the map $B_p$ is defined in equation \ref{eq:beta}).
The values of $a_j$ and $b_j$ will indicate
the $j$-th derivative of $\rho$ at $0$ and $\theta_M$, respectively.
We claim that there exists a unique polynomial
$P$ of degree at most $2p+4$ such that
the following conditions hold:
\[ (\exp \circ P)^{(j)}(0) = a_j, \quad (\exp \circ P)^{(j)}(\theta_M) = b_j, 
\quad j = 0, \ldots, p+1, \]
\[ \int_0^{\theta_M} (\exp \circ P)(\theta) d\theta = 2\pi. \]
This follows from a monotonicity argument analogous
to that used to construct $P_0$ in the case $p = 0$.
Finally, define $\Psi^p(g_0,(\aaa,r))$
to be the potential corresponding to
\[ \rho(\theta) = \exp\left(P(\theta) + 
r(\theta/\theta_M) +
c\,\theta^{p+2}(\theta_M - \theta)^k{p+2}\right) \]
where $c$ is again the unique constant for which
$\int_0^{\theta_M} \rho(\theta) d\theta = 2\pi$.
It is clear that $\Psi^p: G_0 \times H \to H^p(\Ss^1)$
is a diffeomorphism with all the required properties.
\qed

\section{Periodic Sturm-Liouville operators}

For $p \in \ZZ$, $p \ge 0$, and $q \in H^p(\Ss^1)$
we consider the operator $L = L_p(q): H^{p+2}(\Ss^1) \to H^p(\Ss^1)$,
$Lv = -v'' + qv$.
It is easy to verify that $L$ is a Fredholm operator
of index $0$ with kernel of dimension at most $2$.
In particular, the spectrum $\sigma(L)$ is given by
\[ \sigma(L) = \{ \lambda \;|\; \dim\ker(L - \lambda I) > 0 \} \]
and we call $\dim\ker(L - \lambda I)$ the multiplicity
of the eigenvalue $\lambda$.
For $p = 0$ this operator is self-adjoint
and it follows that for all $p \ge 0$ the spectrum of $L$
consists only of real eigenvalues with multiplicity
(geometric equal to algebraic) at most $2$.
We are interested in the geometry of the triple 
$(\Cc_0, \Cc_1, \Cc_2)$ where $\Cc_0 = H^p(\Ss^1)$ and
\[ \Cc_j = \{ q \in H^p(\Ss^1) \;|\; \dim\ker L_p(q) \ge j \}. \]
Recall that $Z(G) = \{\iota^k, k \in \ZZ\}$, the center of $G$,
is the set of vertices of the cones
in $T_{\pm 2}$ (see figure \ref{fig:trlevel}).
The diffeomorphism $\Psi^p_{\Ss^1}$ is the one constructed
in theorem \ref{theo:geomonodromy}.

\begin{theo}
\label{theo:triplesp}
For any $p \in \ZZ$, $p > 0$, $\Psi^p_{\Ss^1}$ is a diffeomorphism
from the triple $(G_0, T_2 \cap G_0, Z(G) \cap (T_2 \cap G_0)) \times \HH$
to $(\Cc_0, \Cc_1, \Cc_2)$.
\end{theo}

{\nobf Proof: }
In a nutshell, potentials whose monodromy is in $T_2$ (resp., its vertices)
belong to $\Cc_1$ (resp., $\Cc_2$).
More precisely, given a potential $q \in H^p(\Ss^1)$,
\[ q \in \Cc_1 
\quad\Longleftrightarrow\quad
\mu(q) \textrm{ has eigenvalue } 1
\quad\Longleftrightarrow\quad
\tr(\mu(q)) = 2
\quad\Longleftrightarrow\quad \]
\[ \quad\Longleftrightarrow\quad 
\mu(q) \in T_2 \cap G_0
\quad\Longleftrightarrow\quad
(\Psi^{\Ss^1}_p)^{-1}(q) \in (T_2 \cap G_0) \times \HH. \]
Also,
\[ q \in \Cc_2
\quad\Longleftrightarrow\quad
\mu(q) = \iota^{2k}, k \in \ZZ, k > 0
\quad\Longleftrightarrow\quad \]
\[ \quad\Longleftrightarrow\quad
\mu(q) \in Z(G) \cap (T_2 \cap G_0)
\quad\Longleftrightarrow\quad
(\Psi^{\Ss^1}_p)^{-1}(q) \in (Z(G) \cap (T_2 \cap G_0)) \times \HH. \]
The result is now obvious.
\qed

In particular, the set $\Cc_1$ of potentials $q \in H^p(\Ss^1)$
with $0$ in the spectrum is a disjoint union of
a (topological) hyperplane $\Psi^p_{\Ss^1}((T_2^0 \cap G_0) \times \HH)$
and countably many cones $\Psi^p_{\Ss^1}((T_2 \cap A_n) \times \HH)$, $n > 0$.
Recall that each cone has two sheets, meeting at a vertex,
a topological subspace of codimension $3$.

Let $q_+ \in H^p(\Ss^1)$ be an almost everywhere strictly positive function
and for $q_0 \in H^p(\Ss^1)$, consider the parametrized straight line
$q_0 - s q_+$, $s \in \RR$.
Standard oscillation theory implies the existence of a sequence
of continuous functions $s_i: H^p(\Ss^1) \to \RR$,
\[ s_0(q_0) < s_1(q_0) \le s_2(q_0) < s_3(q_0) \le s_4(q_0) < \cdots \]
such that $0$ is the $n$-th eigenvalue of the potential
$q_0 - s_n(q_0) q_+$.
In particular, $0$ is the (simple) ground state of $q_0 - s_0(q_0) q_+$.

Combining these two points of view, we have the following result.

\begin{theo}
\label{theo:oscillation}
Each straight line $q_0 - s q_+$, $q_0, q_+ \in H^p(\Ss^1)$,
$q_+$ strictly positive a. e., meets the hyperplane
and each sheet of a cone in $\Cc_1$ exactly once.
More precisely,
\[ q_0 - s_n(q_0) q_+ \in \begin{cases}
\Psi^p_{\Ss^1}((T_2^0 \cap G_0) \times \HH),& n = 0,\\
\Psi^p_{\Ss^1}((T_2 \cap A_{\lceil n/2 \rceil}) \times \HH),& n > 0.
\end{cases} \]
\end{theo}

Thus, the $2n-1$ and $2n$-th eigenvalues of $q_0$ coincide
if and only if the line $q_0 + s$, $s \in \RR$,
passes through the vertex of $\Psi^p_{\Ss^1}((T_2 \cap A_n) \times \HH)$.
Also, the set of potentials $q$ for which $0$ is the double eigenvalue
in positions $2n-1, 2n$ is a (topological) subspace of codimension $3$.

As a final application, we describe the critical set of the nonlinear periodic
Sturm-Liouville operator with quadratic nonlinearity.

\begin{coro}
\label{coro:bullet}
Let $p \ge 2$ and $F: H^p(\Ss^1) \to H^{p-2}(\Ss^1)$ be given by
$F(u) = - u'' + u^2/2$. 
Let $C \subset H^p(\Ss^1)$ be the critical set of $F$.
Then the pair $(H^p(\Ss^1), C)$ is diffeomorphic to
$(G_0, T_2 \cap G_0) \times \HH$.
\end{coro}

{\nobf Proof: }
A simple computation shows that
\[ C = \{ u \in H^p(\Ss^1) \;|\; L_{p-2}(u): H^p \to H^{p-2}
\textrm{ has nontrivial kernel} \}. \]
A standard regularity argument shows that for $u \in H^p \subset H^{p-2}$,
$\ker L_{p}(u) = \ker L_{p-2}(u) \subset H^{p+2}(\Ss^1)$
and therefore
\[ C = \{ u \in H^p(\Ss^1) \;|\; L_{p}(u): H^{p+2} \to H^p
\textrm{ has nontrivial kernel} \} \]
which is $\Cc_1$ in the notation of theorem \ref{theo:triplesp},
completing the proof.
\qed

\section{Other boundary conditions}

The results above extend appropriately to other boundary conditions.
For a real $2 \times 4$ matrix $U$ of rank $2$,
let $H^2_U([0,2\pi]) \subset H^2([0,2\pi])$
be the space of functions $v$ satisfying {\it $U$-boundary conditions}:
\[ U \begin{pmatrix} v(0) & v'(0) & v(2\pi) & v'(2\pi) \end{pmatrix}^\ast
= \begin{pmatrix} 0 & 0 \end{pmatrix}^\ast. \]
In particular, $H^2_{(I\; -I)}([0,2\pi]) = H^2(\Ss^1)$ and
$H^2_{(-I \; -I)}([0,2\pi])$ is the space of antiperiodic functions,
where $I$ is the $2 \times 2$ identity matrix.
We shall not discuss higher orders of differentiability
in this setting.

Two classes of matrices $U$ will be of interest:
\[ U_{\theta_0, \theta_{2\pi}}
= \begin{pmatrix} -\sin\theta_0 & \cos\theta_0 & 0 & 0 \\
0 & 0 & -\sin\theta_{2\pi} & \cos\theta_{2\pi} \end{pmatrix},
\quad U_A = (A \; -I) \]
where $\theta_0 \in [0,\pi)$, $\theta_{2\pi} \in (0,\pi]$
and $A \in GL(2,\RR)$.
Set $L_U: H^2_U([0,2\pi]) \subset H^0([0,2\pi]) \to H^0([0,2\pi])$,
$L_U(v) = -v'' + q v$ where $q \in H^0([0,2\pi])$ is a real potential.
In either case, it is easy to verify that $L_U$ is a Fredholm operator
of index $0$ with kernel of dimension at most $2$.
Indeed, $L_U$ is the composition of the inclusion
$H^2_U([0,2\pi]) \subset H^2([0,2\pi])$,
the invertible map $u \mapsto (-u''+qu, u(0), u'(0))$
from $H^2([0,2\pi])$ to $H^0([0,2\pi]) \times \RR^2$
and the projection onto the first coordinate.
As in the periodic case, we consider the triple
$(\Cc_0(U), \Cc_1(U), \Cc_2(U))$ where
\[ \Cc_j(U) = \{ q \in H^0([0,2\pi]) \;|\; \dim\ker L_U \ge j \}. \]

As is well known, the operator $L_{U_{\theta_0, \theta_{2\pi}}}$
is self-adjoint with spectrum of the form
$\lambda_0 < \lambda_1 < \lambda_2 < \ldots$.
Thus, $\Cc_2(U_{\theta_0, \theta_{2\pi}}) = \emptyset$.
Also, $\Cc_1$ has countably many components, each of them a hyperplane.
Indeed, let $v = (\cos\theta_0, \sin\theta_0) \in \RR^2$ and
\[ C_n = \{g \in G_0 \;|\; \arg(v,gv) = \theta_{2\pi} - \theta_0 + n\pi \}
\subset G_0. \]
Here $\arg(v,gv)$, $v \in \RR^2$, $g \in G$,
denotes the angle between $v$ and $gv$.
More precisely, let $\gamma: [0,1] \to G$, $\gamma(0) = I$, $\gamma(1) = g$;
define a continuous function $\alpha: [0,1] \to \RR$ so that $\alpha(0) = 0$
and $\alpha(t)$ is the angle between $v$ and $\gamma(t)v$;
we define $\arg(v,gv)$ to be $\alpha(1)$.
It is easy to see that the sets of lifted matrices $C_n \subset G_0$
are disjoint topological hyperplanes
and that there exists a diffeomorphism from $G_0$ to $\RR^3$
taking each $C_n$ to the plane $z = n$.
The diffeomorphism $\Psi_{[0,2\pi]}^0$ from $G_0 \times \HH$ to $H^0([0,2\pi])$
takes $C_n \times \HH$ to the component of $\Cc_1$ of potentials $q$
such that $0 = \lambda_n$.
Summing up, $(\Cc_0, \Cc_1)$ is diffeomorphic to $(\RR, \NN) \times \HH$.
Oscillation theory is rather simple:
lines of the form $q_0 - s q_+$, $q_+ > 0$ a. e.,
meet each component of $\Cc_1$ exactly once, transversally.


In the case $U = (A \; -I)$, $L_U$ is self-adjoint
if and only if $A \in SL(2,\RR)$.
The geometry of the triple $(\Cc_0, \Cc_1, \Cc_2)$ is now subtler.
We begin by relating the existence of 
a solution satisfying $U$-boundary conditions
to an algebraic property of $\mu(q)$.

\begin{prop}
\label{prop:triplesAalgebra}
Let $U = (A \; -I)$.
The homogeneous equation \ref{eq:ast}
admits a solution satisfying $U$-boundary conditions
if and only if $\tr(M_1) = a + (\sgn\det A)/a$,
where $a = \sqrt{|\det A|}$ and
$M_1 = a A^{-1} \mu(q) \in SL^\pm(2,\RR)$.
Also, all solutions of the homogeneous equation 
satisfy $U$-boundary conditions if and only if
$A \in SL(2,\RR)$ and $A^{-1} \mu(q) = I$.
\end{prop}

{\nobf Proof: }
Given a potential $q \in H^0([0,2\pi])$,
the following conditions are equivalent:
\begin{itemize}
\item{the homogeneous equation \ref{eq:ast}
admits a solution satisfying $U$-boundary conditions;}
\item{there is a nonzero vector $v \in \RR^2$ such that
$\mu(h) v = A v$;}
\item{$1$ is an eigenvalue of $A^{-1} \mu(h) \in GL(2,\RR)$;}
\item{$a$ is an eigenvalue of $M_1 = a A^{-1} \mu(h)$;}
\item{$\tr(M_1) = a + (\sgn\det A)/a$,
where $M_1 = a A^{-1} \mu(h) \in SL^\pm(2,\RR)$.}
\end{itemize}
This implies the first claim. As to the second claim,
it is clear that all solutions satisfy $A$-boundary conditions
if and only if $A = \mu(q)$.
\qed

\begin{theo}
\label{theo:triplesA}
Let $U = (A \; -I)$.
If $\det A < 0$, then $\Cc_2(U) = \emptyset$ and
$\Cc_1(U)$ is a topological hyperplane.
If $\det A > 0$, $\det A \ne 1$, then $\Cc_2(U) = \emptyset$ and
there exists a diffeomorphism between the pairs
$(\Cc_0(U), \Cc_1(U))$ and
$(G_0, T_{\pm 4} \cap G_0) \times \HH$.
The triples
$(\Cc_0(U), \Cc_1(U), \Cc_2(U))$ and
$(G_0, T_{\pm 2} \cap G_0, Z(G) \cap (T_{\pm 2} \cap G_0)) \times \HH$
are diffeomorphic if $\det(A) = 1$.
\end{theo}

{\nobf Proof: }
From proposition \ref{prop:triplesAalgebra},
the diffeomorphism $(\Psi^0_{[0,2\pi]})^{-1}$ takes the triple
$(\Cc_0(U), \Cc_1(U), \Cc_2(U))$ to
$(G_0, C_1(A), C_2(A)) \times \HH$ where
\[ C_1(A) = \{ M \in G_0 \;|\;
M_1 = a A^{-1} M, \tr(M_1) = a + (\sgn\det A)/a \}; \]
$C_2(A)$ is empty if $A \not\in SL(2,\RR)$, and, if $A \in SL(2,\RR)$, 
\[ C_2(A) = \{ M \in G_0 \;|\; M_1 = a A^{-1} M, M_1 = I \}. \]
It suffices to characterize the triple $(G_0, C_1(A), C_2(A))$
up to diffeomorphism.

Set $B = a A^{-1} = \phi^\pm_L(\theta,\ast,\ast)$ and
define $\beta: G \to G^\pm$, $\beta(M) = BM = M_1$.
We claim that $\beta(G_0) = G^\pm_{\theta}$.
Indeed, for $\det(A) > 0$,
\[ M \in G_0
\quad\Longleftrightarrow\quad
\arg(e_2, Me_2) < 0
\quad\Longleftrightarrow\quad
\arg(Be_2, BMe_2) < 0
\quad\Longleftrightarrow\quad \]
\[ \quad\Longleftrightarrow\quad 
\arg(e_2, Be_2) + \arg(Be_2, M_1e_2) < \arg(e_2,Be_2)
\quad\Longleftrightarrow\quad \]
\[ \quad\Longleftrightarrow\quad
\arg(e_2, M_1e_2) < -\theta
\quad\Longleftrightarrow\quad
M_1 \in G_\theta. \]

Similarly, for $\det(A) < 0$,
\[ M \in G_0
\quad\Longleftrightarrow\quad
\arg(e_2, Me_2) < 0
\quad\Longleftrightarrow\quad
\arg(\tilde RBe_2, \tilde RBMe_2) < 0
\quad\Longleftrightarrow\quad \]
\[ \quad\Longleftrightarrow\quad 
\arg(e_2, \tilde RBe_2) + \arg(\tilde RBe_2, \tilde RM_1e_2) <
\arg(e_2,\tilde RBe_2)
\quad\Longleftrightarrow\quad \]
\[ \quad\Longleftrightarrow\quad
\arg(e_2, \tilde RM_1e_2) < -\theta
\quad\Longleftrightarrow\quad
\tilde RM_1 \in G_\theta
\quad\Longleftrightarrow\quad
M_1 \in G^-_\theta. \]
Also, $\beta(C_1(A)) = \{ M_1 \in G^\pm_\theta \;|\; \tr(M_1) = a + (1/a) \} =
T_{a+(1/a)} \cap G^\pm_\theta$ and
$\beta(C_2(A)) = Z(G) \cap (T_{a+(1/a)} \cap G_\theta)$.
Thus, $\beta$ is a diffeomorphism from the triple $(G_0, C_1(A), C_2(A))$
to the triple
$(G_\theta, T_{a+(1/a)} \cap G_\theta, Z(G) \cap (T_{a+(1/a)} \cap G_\theta))$.
Since $a + (1/a) \ge 2$ with equality exactly when $A \in SL(2,\RR)$,
proposition \ref{prop:fivepairs} finishes the case $\det(A) > 0$.
The case $\det(A) < 0$ follows from proposition \ref{prop:onepair}.
\qed

Oscillation theory for $A \in SL(2,\RR)$ works as in the periodic case:
the straight lines $q_0 - s q_+$ meet the ground hyperplane in $\Cc_1$
(if it exists) exactly once and each cone in $\Cc_1$ twice,
unless the straight line goes through the vertex.
It is not clear how oscillation theory fits in for the cases
$A \not\in SL(2,\RR)$. For instance, for 
\[ A = \begin{pmatrix} 1 & 0 \\ 0 & -1 \end{pmatrix}, \]
$q_0 = 0$ and $q_+ = 1$, the whole line $q_0 - s q_+$
is contained in $\Cc_1$:
all functions $q \in H^0([0,2\pi])$ satisfying
$q(2\pi - t) = q(t)$ belong to $\Cc_1$.

%

\vfil


\bigskip\bigskip\bigbreak

{

\parindent=0pt
\parskip=0pt
\obeylines

Dan Burghelea, Ohio State University, burghele@math.ohio-state.edu
Nicolau C. Saldanha, PUC-Rio and Ohio State University,
nicolau@mat.puc-rio.br; http://www.mat.puc-rio.br/$\sim$nicolau/
Carlos Tomei, PUC-Rio, tomei@mat.puc-rio.br

\smallskip

Department of Mathematics, Ohio State University,
231 West 18th Ave, Columbus, OH 43210-1174, USA

\smallskip

Departamento de Matem\'atica, PUC-Rio
R. Marqu\^es de S. Vicente 225, Rio de Janeiro, RJ 22453-900, Brazil

}

\end{document}